\documentclass[11pt, reqno]{amsart}
\usepackage{amsrefs, amsmath, amssymb, amsthm, graphicx,marvosym}
\usepackage{cancel}

\def\({\left(\begin{array}{cccccc}}
\def\){\end{array}\right)}

\def\ds{\displaystyle}

\def\({\left(\begin{array}{cccccc}}
\def\){\end{array}\right)}

\def\bega{\begin{array}}
\def\enda{\end{array}}

\newcommand{\ba}{\overset{\raisebox{0pt}[0pt][0pt]{\text{\raisebox{-.5ex}{\scriptsize$\leftharpoonup$}}}}}
\newcommand{\fa}{\overset{\raisebox{0pt}[0pt][0pt]{\text{\raisebox{-.5ex}{\scriptsize$\rightharpoonup$}}}}}

\newcommand{\bel}{\begin{equation}\label}
\newcommand{\beq}{\begin{equation}}
\newcommand{\eeq}{\end{equation}}
\newcommand{\bes}{\begin{eqnarray}}
\newcommand{\ees}{\end{eqnarray}}
\newcommand{\besnn}{\begin{eqnarray*}}
\newcommand{\eesnn}{\end{eqnarray*}}
\newcommand{\si}{\ensuremath{\sigma}}

\newcommand{\C}{\mathcal C}
\newcommand{\Ord}{\mathcal O}

\newcommand{\ve}{\varepsilon}

\newcommand{\bp}{\begin{proof}}
\newcommand{\ep}{\end{proof}}

\newtheorem{theorem}{Theorem}[section]

\newtheorem{lemma}[theorem]{Lemma}

\newtheorem{remark}[theorem]{Remark}

\numberwithin{equation}{section}

\begin{document}

\title{No BV bounds for approximate solutions to p-system with general pressure law}

\author{Alberto Bressan}
\address{Department of Mathematics, Penn State University,
University Park, PA.~16802, USA\newline ({\tt bressan@math.psu.edu}).}
\author{Geng Chen}
\address{School of Mathematics,
Georgia Institute of Technology, Atlanta, GA 30332 USA \newline({\tt gchen73@math.gatech.edu}).}
\author{Qingtian Zhang}
\address{Department of Mathematics, Penn State University,
University Park, PA.~16802, USA \newline({\tt zhang\_q@math.psu.edu}).}
\author{Shengguo Zhu}
\address{ Department of Mathematics, Shanghai Jiao Tong University, Shanghai 200240, P.R.China
\newline({\tt zhushengguo@sjtu.edu.cn}).
}

\date{\today}

\begin{abstract}
For the p-system with large BV initial data,
an assumption introduced in \cite{bak} by Bakhvalov guarantees the global existence of
entropy weak solutions with uniformly bounded total variation. The present paper
provides a partial converse to this result.
Whenever Bakhvalov's condition  does not hold, we show that
there exist  front tracking approximate
solutions, with uniformly positive density, whose total
variation becomes arbitrarily large.   The construction extends the arguments
in \cite{BCZ} to a general class of pressure laws.
\end{abstract}

\maketitle


\section{Introduction}
A satisfactory existence-uniqueness theory is now available for hyperbolic systems of conservation laws in one space dimension with small total variation 
\cites{glimm, bressan, bly}.
A major remaining open problem is whether the total variation remains uniformly bound or can blow up in finite time for large BV initial data. Up to now, only  few systems of hyperbolic conservation laws are known, where uniform BV estimates hold for solutions
with large data \cites{nishida, tem83}.
On the other hand, examples with finite time blowup have been constructed
in \cite{bj, jenssen}. However, these systems do not come from physical models and do not admit a strictly convex entropy.

In this paper, we focus on the p-system with general pressure law modeling barotropic gas dynamics.
	\beq\label{ps}\left\{
		\begin{array}{rcl}
		u_t+p(v)_x &=&0\,,\\
		v_t-u_x &=& 0\,,
		\end{array}\right.
	\eeq
where $v=1/\rho>0$ is the specific volume, $\rho>0$ is the density and $u$ is the velocity of the gas. The pressure $p(v)$
is a smooth function of $v$ satisfying
\beq\label{pga}
p_v<0\quad \text{and}\quad p_{vv}>0\,.
\eeq

In \cite{nishida}, Nishida proved the global BV existence to \eqref{ps} with large
initial data,
for $\gamma$-law pressure  $p=v^{-\gamma}$ with $\gamma=1$.
On the other hand, in the case
$\gamma=3$, various front tracking approximate solutions were recently  constructed in
\cite{BCZ},  exhibiting finite time blowup of the BV norm.

For the p-system with general pressure law,  in \cite{bak}, Bakhvalov extended the global BV existence result for isothermal gas dynamics in \cite{nishida} to any pressure law
$p(v)$ satisfying the \emph{Bakhvalov's condition}
\beq\label{bak_con}
3p_{vv}^2~\leq~ 2p_{v}p_{vvv}\quad  \text{for\ all}\quad v> 0\,.
\eeq
In particular, for $\gamma$-law pressure  $p=v^{-\gamma}$ with $\gamma>0$, Bakhvalov's condition holds if and only if $\gamma\in(0,1]$.
In \cite{bak}, more general $2\times2$ systems of conservation laws are also considered.

\begin{figure}[htbp]
\centering
  \includegraphics[scale=0.4]{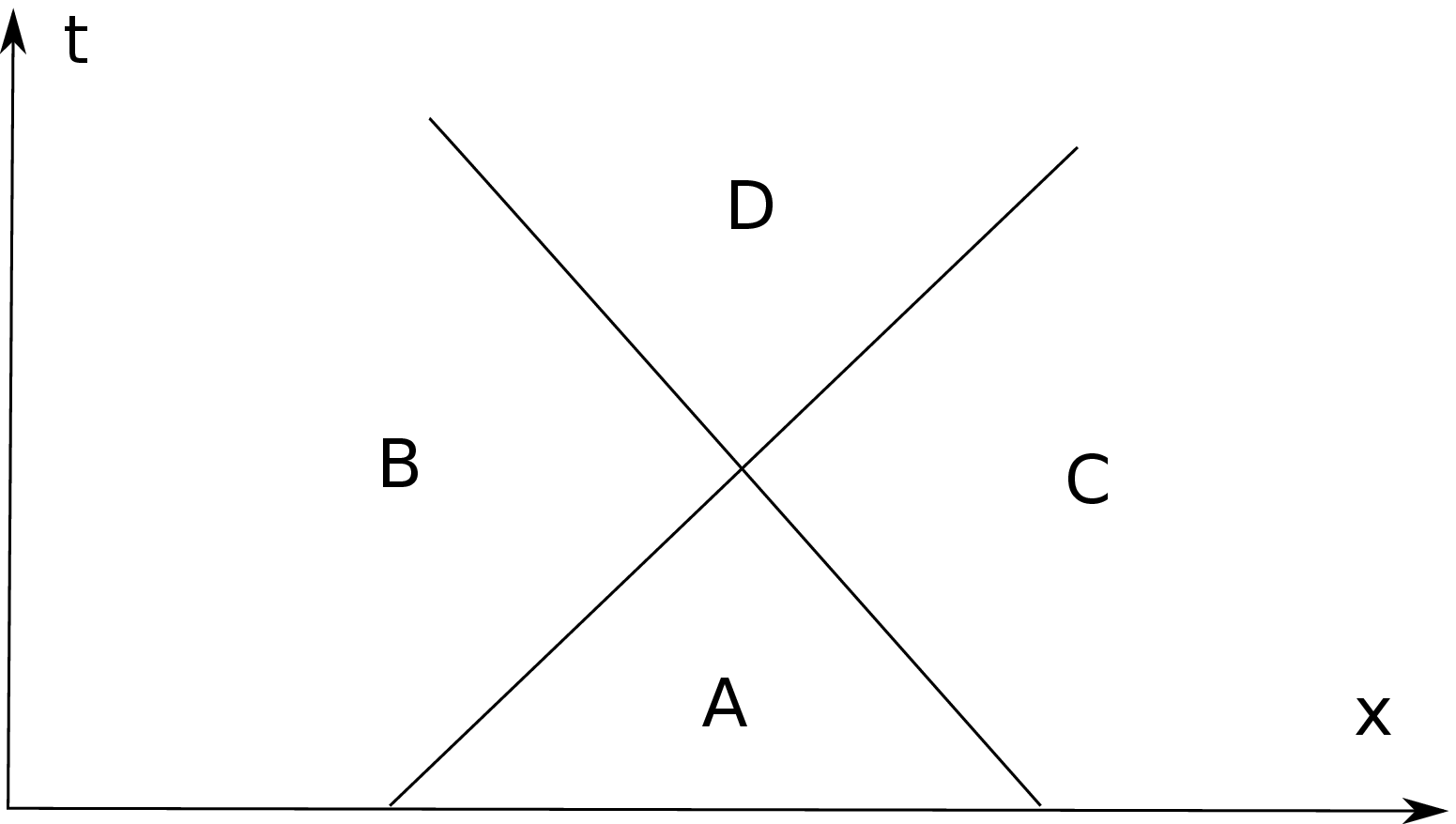}\hfill
  \includegraphics[scale=0.8]{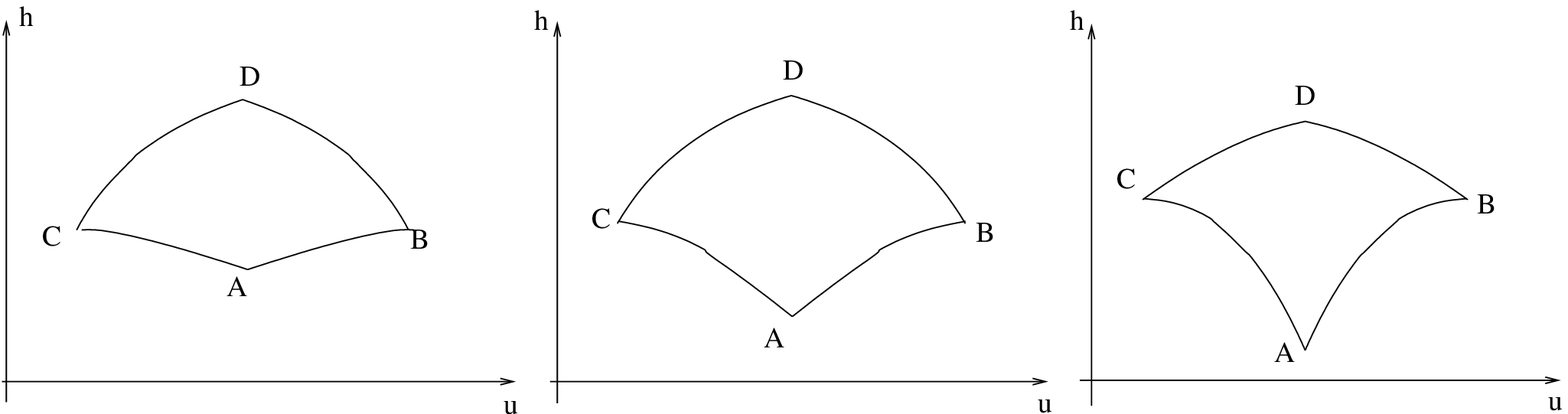}
    \caption{\small The upper figure shows two interacting shocks in the $x$-$t$ plane.
    For this interaction, the
 lower figures show three different cases, in the $(u,h)$-plane
Left: the strength is amplified after the crossing, and Bakhvalov's condition (\ref{bak_con}) is not satisfied.
    Middle: the strength is same after crossing.
     Right: the strength is reduced after crossing.
     In the middle and the right pictures, Bakhvalov's condition (\ref{bak_con}) is satisfied.
     For  $\gamma$-law pressure with $\gamma>0$, the left figure corresponds to
      $\gamma>1$, the middle to $\gamma=  1$, and the right to $0<\gamma<1$.
     }
\label{f:ba}
\end{figure}

We observe that Bakhvalov's condition
 determines whether the strength of a shock increases or decreases
by crossing a shock of the
opposite family, as shown in Figure \ref{f:ba}.
The shock strength is here measured by the change of $h(v)$ across the shock,
where
\bel{defh}
h(v)~\doteq~\int_{v}^1 \sqrt{-{p_v}}\,dv
\eeq
is the density part in the Riemann invariants
\[
s~\doteq~u+h(v),\qquad r~\doteq~u-h(v).
\]

In the present paper we extend the blowup examples in \cite{BCZ}  to the case where
the pressure  violates \eqref{bak_con}.
Together with \cite{bak}, this indicates that Bakhvalov's condition (\ref{bak_con})
is necessary and sufficient for the BV stability  of the front tracking scheme.
More precisely, the following result will be proved.
\vspace{.2cm}

{\bf Theorem.}    {\it Assume that the pressure $p(\cdot)$
satisfies (\ref{pga}) for every $v>0$ but violates \eqref{bak_con}
for some $v>0$.  Then there exists a front tracking approximate solution
where the density remains uniformly positive while the  total strength of waves
approaches infinity as $t\to\infty$.
At each wave-front interaction the strengths of outgoing waves
are the same as in the exact solution. The
only errors introduced by the front tracking approximation
are in the speeds of the wave fronts.}

\vspace{.2cm}

By suitably modifying
the construction given  in the last section of \cite{BCZ}, we expect that one could also
construct an example of front tracking approximation where the  BV norm
blows up in with finite time. The main ideas leading to the blow-up example can be explained as follows.\vspace{.2cm}

\begin{figure}[htbp]
\centering
  \includegraphics[scale=0.5]{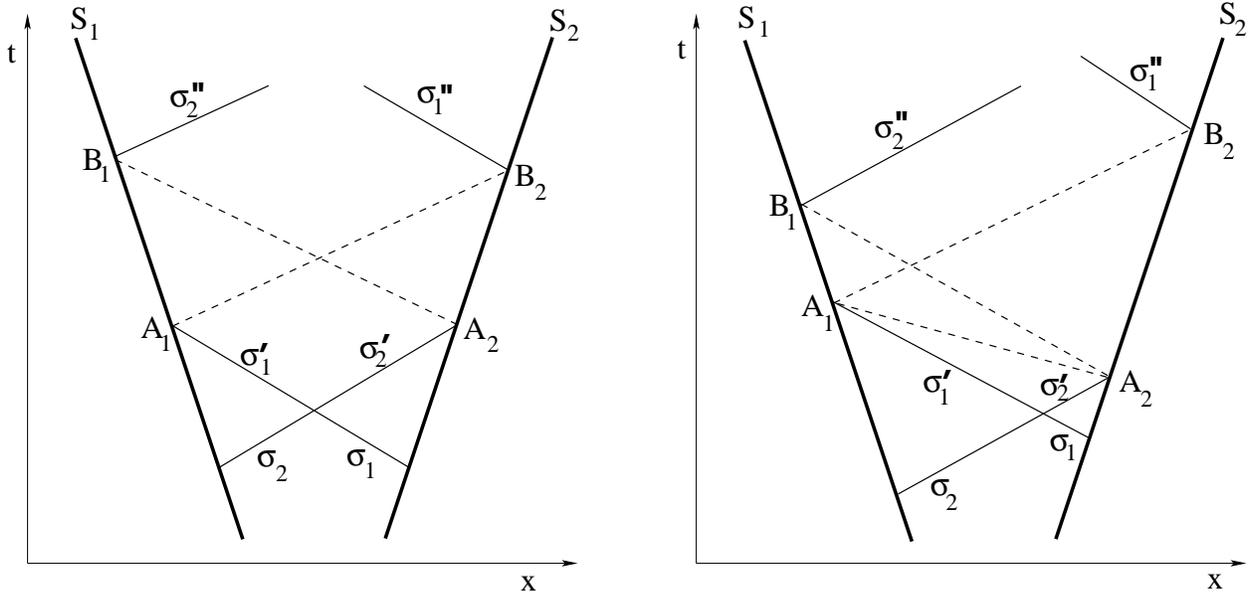}
    \caption{\small Solid lines denote shocks
    while dotted lines denote rarefaction fronts.
     Left: a symmetric periodic interaction pattern, where two small
    fronts bounce back and forth between two large shocks.  In the region between the two
    large shocks, the solution takes values in the region where  Bakhvalov's condition fails.   Right: an asymmetric interaction pattern.   Here part of the rarefaction originating from
    $A_2$ is canceled at $A_1$ by merging with a shock of the same family.   }
\label{f:hyp102}
\end{figure}

If Bakhvalov's condition (\ref{bak_con}) fails for $v$ in a neighborhood of
$v_0$, one can  construct two small approaching shocks such that
\begin{itemize}
\item[(i)] their left and right states remain in the region where Bakhvalov's condition fails,
and hence
\item[(ii)] calling  $\sigma_1,\sigma_2$ their sizes before the interaction and
$\sigma'_1,\sigma'_2$ their sizes after the interaction, one has
\bel{amplif}\sigma_1~=~\sigma_2~<~\sigma_1'~=~\sigma_2'\,.\eeq
\end{itemize}
Next, assume that these small fronts bounce back and forth between two
very large shocks (Fig.~\ref{f:hyp102}, left).   After a first reflection at the points
$A_1,A_2$, two rarefactions are
created.   When these rarefaction impinge again on the large shocks at $B_1,B_2$,
they generate
two new shocks.   Every time a front is reflected by a large shock, the outgoing wave is strictly smaller than the incoming one.    However, is the shocks $S_1,S_2$ are
very large, the strengths of incoming and reflected fronts are almost the same.
Thanks to (\ref{amplif}), by a suitable choice of the shock strengths, we can achieve
\bel{perpat1}\sigma_1~=~\sigma_2~=~\sigma_1''~=~\sigma_2''\,.\eeq
Hence the interaction pattern can be iterated in time.

If we further increase the strengths of the shocks $S_1,S_2$, in (\ref{perpat1})
we would have
\bel{perpat2}|\sigma_1|~=~|\sigma_2|~<~|\sigma_1''|~=~|\sigma_2''|\,.\eeq
To achieve again a periodic pattern, one needs to cancel part
of the rarefaction emerging at $A_2$.   As shown in Fig.~\ref{f:hyp102}, right,
this can be done by merging it with a shock
of the same family, at the interaction point $A_1$.
In the end, this yields an asymmetric, periodic
interaction pattern where
$$\sigma_1''~=~\sigma_1\,,\qquad\qquad \sigma_2''~=~\sigma_2\,.$$

\begin{figure}[htbp]
\centering
  \includegraphics[scale=0.5]{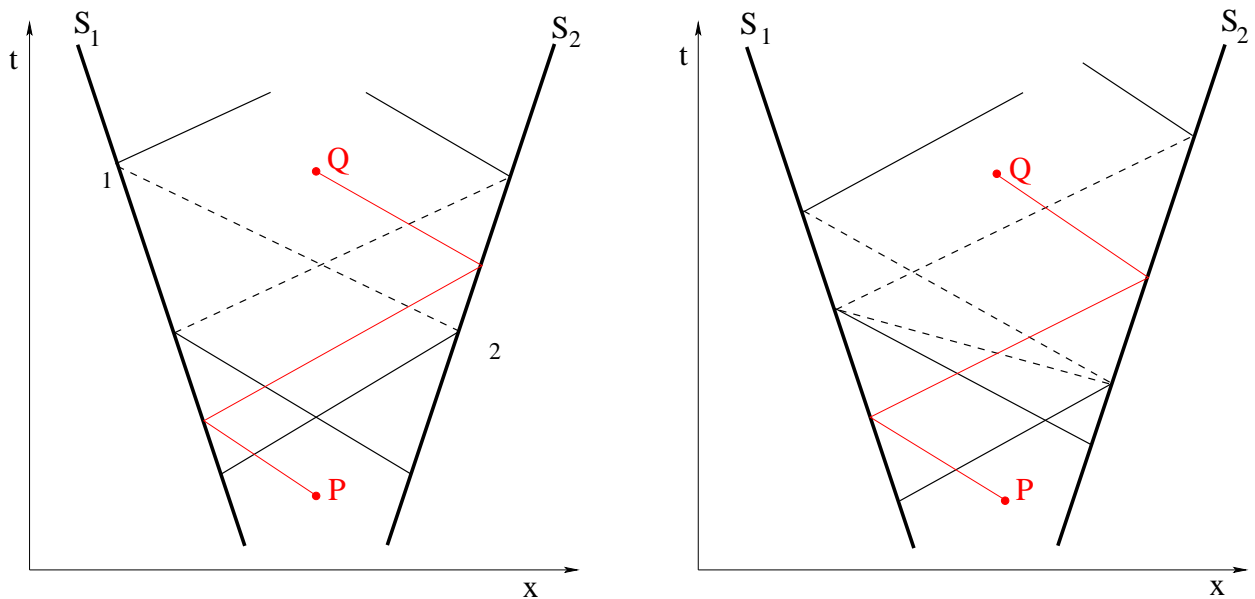}
    \caption{\small If a small front is added on top of the interaction patterns in
    Fig.~\ref{f:hyp102},  after a complete set of interactions the strength of this
    front is  (i) almost the same,  in case of the symmetric pattern on the left, and
    (ii) strictly larger, in case of the asymmetric pattern on the right. }
\label{f:hyp103}
\end{figure}

Next, on top of  these periodic patterns we add an infinitesimally small wave front
(a compression or a rarefaction), as in  Fig.~\ref{f:hyp103}.
If this additional front is initially located at $P$ and has strength $\ve$, after a complete set of interactions
we show that
\begin{itemize}
\item[(i)]
 For the symmetric interaction pattern the strength of the small front
at $Q$ is  $\ve' = \ve + o(\ve)$.
\item[(ii)]
 For the asymmetric interaction pattern the strength of the small front
at $Q$ is  $\ve' = \kappa \ve + o(\ve)$, for some $\kappa>1$.
\end{itemize}
As this cycle of interactions is repeated over and over, 
the infinitesimal front is enlarged by an arbitrarily large factor.

Finally, as in \cite{BCZ},  we  replace this infinitesimally small 
front with a train of countably many pairs
of rarefaction-compression fronts having sizes $\pm 2^{-k}\ve$, with $k=1,2,\ldots$.
This yields a front-tracking approximate solution satisfying the properties
stated in the Theorem (see Fig.~\ref{f:hyp100}).

\vspace{.2cm}
{\bf Remark.}
We emphasize that our result does NOT imply that the total variation of
entropy weak solutions to the p-system can become arbitrarily large.
Rather, it shows that
front tracking approximations can be unstable in the BV norm, whenever
Bakhvalov's condition is violated. The present construction also shows that
for large initial data,
uniform a priori bounds on the total variation cannot be proved simply
by estimating the wave strengths at each interaction.
As remarked in \cite{BCZ}, to establish such BV bounds (if they do indeed hold)
it will be essential to use also  the decay of rarefaction waves, due to genuine nonlinearity.
\vspace{.2cm}

The paper is organized as follows. In Sections~2 and 3 we
study the wave curves and calculate wave interactions.
In Section 4 we first construct a front tracking approximate solution with a
periodic interaction pattern. Then, by suitably perturbing this periodic pattern,
we give examples of front tracking
approximate solutions  where the BV-norm blows up as $t\to\infty$.

\section{Wave curves}
In this section, we introduce basic notation and review the  rarefaction, compression, and shock curves for \eqref{ps}. We omit some standard calculations
for wave curves and refer the reader to Chapter 17 in \cite{smoller} for details.

In  Lagrangian coordinates, the wave speed for \eqref{ps} is
\[
c=\sqrt{-{p}'(v)}.
\]
Integrating the eigenvectors of  \eqref{ps}, one obtains the Riemann invariants $s$ and $r$:
\beq\label{s_r_eq}
s~\doteq~u+h \qquad r~\doteq~u-h\,,
\eeq
with
\beq\label{xi_eq}
 h\equiv h(v)~\doteq~\int_{v}^1 \sqrt{-{p_v}}\,dv\,.
\eeq
In the case of a smooth solution, $s$ and $r$ satisfy
\[
s_t+cs_x~=~0\,,\qquad \qquad r_t-cr_x~=~0\,.
\]
This yields the curves for the Rarefaction and Compression simple waves.

For a shock wave, the Rankine-Hugoniot jump conditions take the form
\bes
\si[u]&=&[p(v)]\,,\\
\si[v]&=&-[u]\,.
\ees
where $[u]=u_r-u_l$, etc$\ldots$, and the subscripts $l$ and $r$ denote the left and right states on the shock wave, respectively. Together with the
Lax entropy condition, this uniquely determines the shock curves.

The following table summarizes the equations for rarefaction, compression and shock curves.
We refer the reader to Chapter 17 in \cite{smoller} for detailed calculations.
We use $(\bar u,\bar v)$ and $(u,v)$ to denote the left and right states across the
wave, respectively. Moreover, we use $\fa R$, $\ba R$,  $\fa C$, $\ba C$ $\fa S$,
and $\ba S$ to denote the forward or backward (or second or first) rarefaction, compression and shock waves respectively.
\beq\label{curve}\left.
\begin{array}{rcl}
	\ba R:\ &u-\bar u=h(\bar{v})-h(v),&v>\bar v\vspace{.1cm}\\
	\fa R:\ &u-\bar u=h({v})-h(\bar v),&v<\bar v\vspace{.1cm}\\
	\ba C:\ &u-\bar u=h(\bar{v})-h(v),&v<\bar v\vspace{.1cm}\\
	\fa C:\ &u-\bar u=h({v})-h(\bar v),&v>\bar v\vspace{.1cm}\\
	\ba S:\ &u-\bar u=-\sqrt{(v-\bar v)(p(\bar v)-p(v))},&v<\bar v\vspace{.1cm}\\
	\fa S:\ &u-\bar u=-\sqrt{(v-\bar v)(p(\bar v)-p(v))},&v>\bar v\,.
\end{array}\right.
\eeq
We recall that  the  combined shock-rarefaction curves have  $\C^2$ regularity
\cite{bressan, smoller}.
%
\section{Wave interactions}
In this section, we calculate the head-on interactions between two shocks and
between a shock and a rarefaction, respectively.

\subsection{Preliminaries}\label{sub_3.1}
Consider the function
\beq\label{a_def}
a~=~ a(v,\bar v)~\doteq~h(\bar v)-h(v)~=~\int^v_{\bar v}\sqrt{-p_v}dv\,.
\eeq
Since $a_v(v,\bar v)>0$, one can recover $v$ as a function of $a$ and $\bar v$,
say, $v=v(a,\bar v)$. We also introduce the function
\beq\label{RH1}
F(a,\bar v)~\doteq~
\sqrt{\Big(v(a,\bar v)-\bar v\Big)\Big(p\big(\bar v\big)-p\big(v(a,\bar v)\big)\Big)}\,.
\eeq
We now compute the Taylor expansion of $p(v)$, for $v$ near $\bar v$.
In turn, this can be used to  calculate the
Taylor expansion of $F(a,\bar v)$.
\beq\label{p_tal}
p(v)~=~p(\bar v)+p'(\bar v)(v-\bar v)+\frac{1}{2}p''(\bar v)(v-\bar v)^2+\frac{1}{6}p'''(\bar v)(v-\bar v)^3
+\frac{1}{24}p^{(4)}(\bar v)(v-\bar v)^4+o(v-\bar v)^4\,.
\eeq
Using (\ref{a_def}) and considering $v=v(a,\bar v)$, we compute
\bes\label{v_tal}
v-\bar v&=&v(a,\bar v)-v(0,\bar v)\\
	    &=&(-p'(\bar v))^{-\frac{1}{2}}a+\frac{1}{4}(-p'(\bar v))^{-2}p''(\bar v)a^2
	    +\frac{1}{6}\Big[(p''(\bar v))^{2}-\frac{1}{2}p'(\bar v) p'''(\bar v)\Big](-p'(\bar v))^{-\frac{7}{2}}a^3+\Ord (a^4)\,.\nonumber
\ees
Using \eqref{p_tal} and \eqref{v_tal}, we obtain
\bes\label{F_tal}
F(a,\bar v)&=&|v-\bar v|\cdot\sqrt{\frac{p(\bar v)-p(v)}{v-\bar v}}\\
	    &=&\text{sign}(v-\bar v)\cdot a\left\{1+J_1(\bar v)a^2+J_2(\bar v) a^3\right\}+o(a^4),
\nonumber
\ees
where
\beq\label{J1_def}
J_1~\doteq~\frac{1}{96}(-p')^{-3}(p'')^{2}\,,
\eeq
\beq\label{J2_def}
J_2~\doteq~\frac{1}{32}p''\big( \frac{1}{2}(p'')^2-\frac{1}{3}p'''p'\big)\,.
\eeq
%
\subsection{Head-on wave interactions}
\begin{figure}[htbp]
\centering
  \includegraphics[scale=0.5]{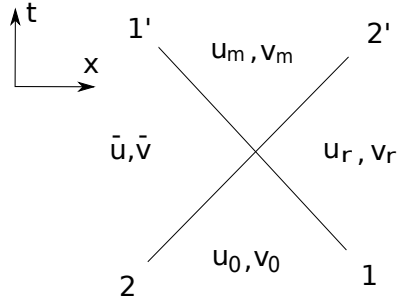}
    \caption{ Head-on interactions}
    \label{f:head-on}
    \end{figure}
    In this section, we consider interactions  between two opposite waves, as shown in
    Figure~\ref{f:head-on}, where the incoming waves can be rarefaction, shock, or
    compression waves. The wave does not change its type after crossing a wave
    of the opposite family.

    We use subscripts $1$, $2$, $1'$ and $2'$ to denote the incoming and outgoing
    waves of the first and second family, respectively.
    And we denote the $(u,v)$ states between these waves according to Figure \ref{f:head-on}.
     For any wave-front, we denote by use $a=h_{\rm left}-h_{\rm right}$
     the difference between the values of $h$ at the left and right states of the
     front. For example, referring to Figure~\ref{f:head-on}, one has
     \beq\label{a_SS}
     a_2=\bar h-h_0,\qquad a_1=h_0-h_r,\qquad  a_{1'}=\bar h-h_m\,,\qquad
     a_{2'}=h_m-h_r\,.
     \eeq

\paragraph{\bf Shock-Shock interaction}
In this part, we consider the weak shock-shock interaction. For simplicity, we only consider the case when two shocks have same strength, i.e. $\bar v=v_r$, hence $\bar h=h_r$.
Using (\ref{curve}) and \eqref{a_SS} with $\bar h=h_r$, one obtains 
\beq\label{a12}
\left\{\bega{rl}a_2&=~\bar h-h_0~=~h_r-h_0~=~-a_1~>~0\,,
\cr\cr
a_{2'}&=~h_m-h_r~=~h_m-\bar h~ =~-a_{1'}~>~0.\enda\right.\eeq
Then by (\ref{curve}) and (\ref{RH1}), we have
\[F(a_2,\bar v)~=~F(a'_1,\bar v)\]
which yields
\beq\label{ss}
a_{2'}~=~-a_{1'}~=~a_2\Big(1+2J_2(\bar v)a_2^3\Big)+o(a_2^4),
\eeq
where $J_2$ was defined at (\ref{J2_def}).
%
\subsection{Rarefaction-Shock or Compression-Shock interaction}
We now consider the interaction between a backward rarefaction and a forward shock.
By \eqref{curve}, we know $a_1$, $a_2$, $a_{1'}$, and $a_{2'}$ are all negative.
Traversing the waves before and after interaction yields
\[
-a_{1'}+F\Big(a_{1'},\bar v\Big)~=~-a_1+F\Big(a_{1'},v_0(a_2,\bar v)\Big).
\]
By the equation \eqref{F_tal}, we thus have
\[
a_{1'}~=~a_1+J_2 a_2a_1^3 +o(a_2a_1^3).
\]
Hence, by
\[
a_1+a_2~=~a_{1'}+a_{2'}
\]
we obtain
\beq\label{sr}
|a'_2|~=~|a_2|(1+J_2|a_1|^3) +o(|a_2a_1^3|)
\eeq

By an entirely similar calculation, we have same estimate for the interaction between
a backward compression and a forward shock.
By symmetry, a similar estimate holds for the interaction between
a forward compression and a backward shock.

%
\section{Front tracking approximations with unbounded BV norm}
In this section, we construct a front tracking approximate solution whose BV norm
tends to infinity as $t\to\infty$. We assume that
the Bakhvalov condition \eqref{bak_con} fails at some $v>0$.
Hence, by continuity and by \eqref{pga} there exists some interval $(v_L,v_U)$, in which
\beq\label{blowup}
J_2~\doteq~\frac{1}{32}p''\big(
3p_{vv}^2- 2p_{v}p_{vvv}\big)~>~0\qquad~~~\hbox{for all}~~v\in (v_L,v_U).
\eeq

\subsection{Front tracking approximations with a periodic interaction pattern.}
\begin{figure}[htbp]
\centering
  \includegraphics[scale=0.4]{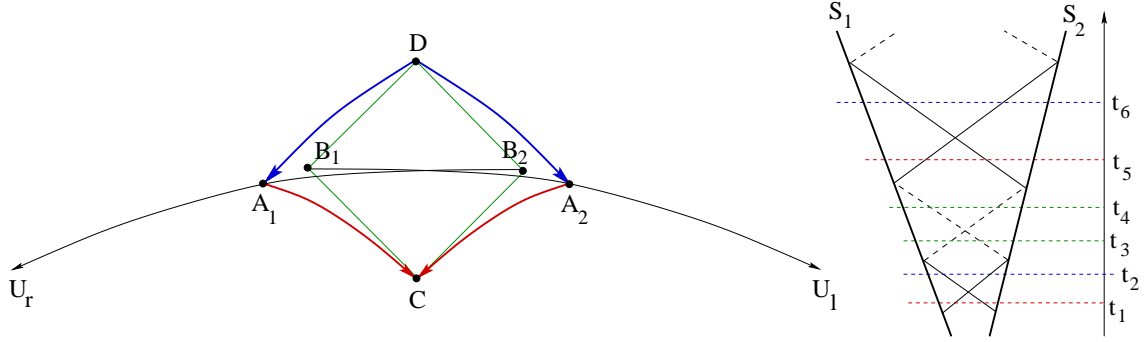}
    \caption{{\small A periodic interaction pattern. The left picture is on $h$-$u$ plane.
    The right picture is on $t$-$x$ plane}}
\label{f:hyp59}
\end{figure}

Following \cite{BCZ}, we  first construct a symmetric interaction pattern containing four wave fronts, as shown in Fig.~\ref{f:hyp59}.
This pattern is symmetric, because two boundary shocks $S_1$ and $S_2$ (and also the inner shocks $A_1C$ and $A_2C$) are chosen to have the same strength measured by the difference in $h$ between two sides of each shock.  We choose the strengths of the two
large shocks $S_1, S_2$ and of the two intermediate waves in such a way that,
after a whole round of interactions, these strengths are the same as at the initial time.  Working in the $(u,h)$ plane, this is achieved as follows.

\begin{itemize}
\item[(i)]  Choose states $A_1$, $A_2$, $B_1$, $B_2$, $C$ and $D$ such that $v\in(v_L,v_U)$ at these states. Hence
		\eqref{blowup} is satisfied inside and on a neighborhood of the diamond
		with vertices $A_1$, $C$, $A_2$, $D$.
\item[(ii)]  Construct two shocks: the 1-shock
$A_1C$ and the 2-shock $A_2C$, approaching each other.
\item[(iii)] Determine the two outgoing shocks $DA_1$ and $DA_2$, resulting
from the crossing of the above two shocks.
\item[(iv)] Construct a rectangle  having two opposite vertices
at $C$ and $D$.   Call $B_1$, $B_2$ the remaining two vertices.
\item[(v)] Finally, the state  $U_l$ is chosen so that the two points $B_1$ and $A_2$
are on the same 1-shock curve with left state $U_l$.
Symmetrically,  $U_r$ is chosen so that the two points $B_2$ and $A_1$
are on the same 2-shock curve with  right state $U_r$.
\end{itemize}
We observe that, by \eqref{ss} and \eqref{blowup}, the $h$-component of the states
 $B_1$ and $B_2$ is larger than the $h$-component of $A_1$ and $A_2$.

\begin{figure}[htbp]
\centering
  \includegraphics[scale=0.5]{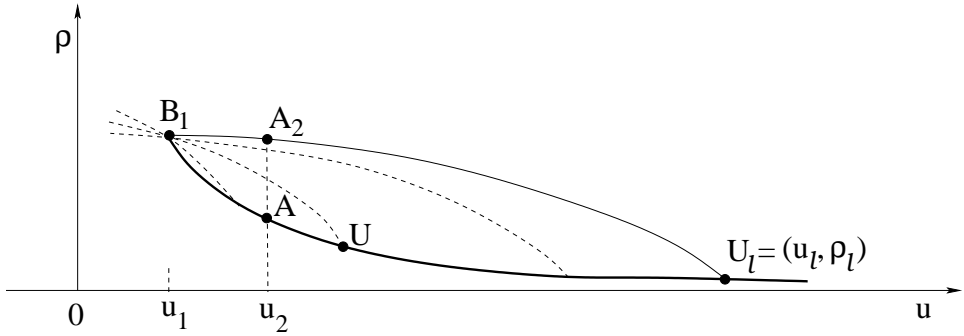}
    \caption{{\small By moving the point $U$ along the 1-shock
    curve with right state $B_1$, we eventually
   reach a left state $U_l$     such that the 1-shock curve through $U_l$
   contains $A_2$ as well. }}
\label{f:hyp62}
\end{figure}

The existence of  states $U_l$, $U_r$ satisfying (v) is now proved in the following lemma,  illustrated in Fig.~\ref{f:hyp62}.

\begin{lemma}
In the $(u,h)$-plane, consider two points $B_1=(u_1,h_1)$
and $A_2=(u_2, h_2)$. Assume that
\begin{itemize}
\item[(i)]~~$u_1<u_2$, and $h_1>h_2$.
\item[(ii)] Calling $A=(u_2, h_2^*)$ the point on the 1-shock curve with right state
 $B_1$
with the same $u$-component as $A_2$, one has $h_2^*<h_2$.
\end{itemize}
Then there exists
a unique
$U_l= (u_l,h_l)$, with $0<h_l<h_2$, such that both $B_1$ and $A_2$ lie on the 1-shock curve
with left state state $U_l$.
\end{lemma}
\begin{remark}{\rm
Condition (ii) clearly holds when the interaction diamond $A_1$-$C$-$A_2$-$D$ is small enough, i.e. the interactions inside the diamond are all weak.}
\end{remark}
\begin{proof}

We shall use (\ref{RH1}) with
$(u_l, v_l)$ while $(u, v) = (u_1, v_1) $ or $ (u_2, v_2) $.
To prove the lemma we need to find $(u_l, \rho_l)$
such that
\beq\label{url}
u_l-u_1~=~\sqrt{\big(p(v_1)-p(v_l)\big)(v_l-v_1)}\,,
\qquad
u_l-u_2~=~\sqrt{\big(p(v_2)-p(v_l)\big)(v_l-v_2)}\,.
\eeq
This will be achieved if we can find $v_l$ such that
\beq\label{url2}
u_2-u_1~=~G(v_l)\,,
\eeq
where $G$ is the function defined as
$$
G(v)~\doteq~\sqrt{\big(p(v_1)-p(v)\big)(v-v_1)}-\sqrt{\big(p(v_2)-p(v)\big)(v-v_2)}\,.
$$
The assumption (ii) implies
$$G(v_2) ~=~ \sqrt{\big(p(v_1)-p(v_2)\big)(v_2-v_1)}~<~u_2-u_1\,.$$
Moreover, a direct computation shows
$$\lim_{v_l\to \infty} G(v_l) ~=~+\infty.$$
Finally, for any $v_1<v_2<v$, we have
\[
\frac{\partial}{\partial v} G(v)~ =~\frac{-p'(v)(v-v_1)+p(v_1)-p(v)}{2 \sqrt{\big(p(v_1)-p(v)\big)(v-v_1)}}-
\frac{-p'(v)(v-v_2)+p(v_2)-p(v)}{2 \sqrt{\big(p(v_2)-p(v)\big)(v-v_2)}}~>~0.
\]
Indeed,  since $p''(v)>0$, one has
\[
\frac{\partial}{\partial a}\Big(\frac{-p'(v)(v-a)+p(a)-p(v)}{2 \sqrt{\big(p(a)-p(v)\big)(v-a)}}\Big)
~=~\frac{\Big(p'(a)-\frac{p(v)-p(a)}{v-a}\Big)\Big(p'(v)-\frac{p(v)-p(a)}{v-a}\Big)}{\big(p(a)-p(v)\big)^3(v-a)^5}~<~0
\]
for any $a<v$.
Since $v_l\geq v_2$, there exists a unique value of $v_l$ such that $G(v_l) = u_2- u_1$.
\end{proof}

\begin{figure}[htbp]
\centering
 \includegraphics[scale=0.5]{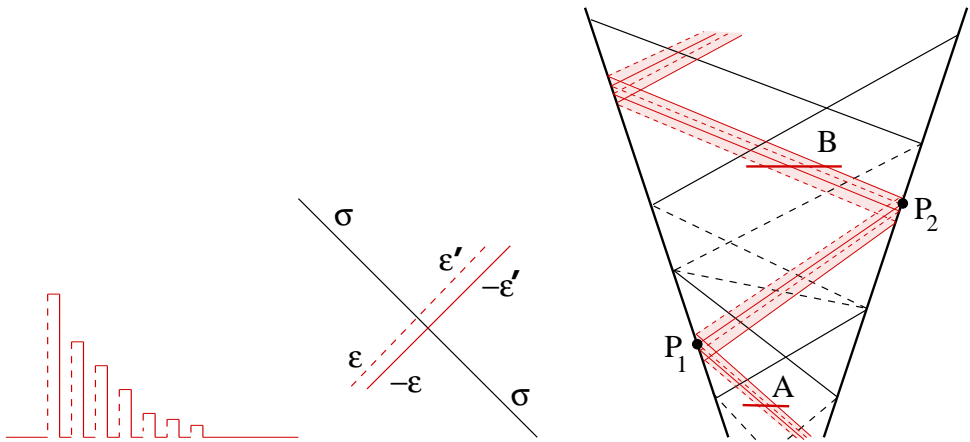}
    \caption{{\small A periodic pattern that amplifies a train of small wave fronts.}}
\label{f:hyp100}
\end{figure}
\bigskip
%

\subsection{An example with unbounded BV-norm}
Next,  as shown in Fig.~\ref{f:hyp100}, on top of the periodic pattern constructed in
Fig.~\ref{f:hyp59},
we add a train  of countably many pairs of rarefaction and compression waves.
The $k$-th pair of waves have sizes
$\pm 2^{-k}\ve$.   Notice that if a front of arbitrary size $\sigma$  crosses
a rarefaction and then a compression wave  of exactly opposite sizes, after the two crossings the size of the front is still $\sigma$, exactly as before  (Fig.~\ref{f:hyp100}, center).
As a result, the interaction pattern of four large fronts retains its periodicity.

Note that in Fig.~\ref{f:hyp100}, we perturb the symmetric periodic pattern in
Fig.~\ref{f:hyp59} to an asymmetric periodic pattern  by splitting some reflecting rarefaction wave into two pieces. The detail of this perturbation will be discussed later.
We recall that the strength of a wave is always defined as
\[
|a(v_{\rm left},v_{\rm right})|~=~|h_{\rm right}-h_{\rm left}|\,,
\]
where the subscripts  denote the left and right states across the wave-front, respectively.

We always assume that each front in the train of small waves has strength
$\leq \ve$. Indeed, we can always
perform a partial cancellation of the compression-rarefaction pair so that
both fronts have strength $ \leq \ve$.  We choose
$\ve>0$  small enough so that all states between two boundary shocks satisfy $v\in(v_L,v_U)$.

We consider the amplification of total wave strength of these alternating waves.
To fix the ideas, consider a 1-rarefaction or compression
of strength $\ve_A>0$, located at $A$.
 Within a time period, this front will
\begin{itemize}
\item[i.] Cross the intermediate 2-shock.
\item[ii.] Interact with the large 1-shock at $P_1$ producing a 2-compression.
\item[iii.] Cross the intermediate 1-shock.
\item[iv.] Cross the intermediate 1-rarefaction.
\item[v.] Interact with the large 2-shock at $P_2$ producing a 2-rarefaction.
\item[vi.] Cross the intermediate 2-rarefaction.
\end{itemize}

Indeed, when a small wave of strength $\ve^-$ crosses a shock of the opposite family
of strength $s$, by (\ref{sr}) the strength of the outgoing front is
\bel{amp1}\ve^+ ~=~\left(1+ J_2s^3+o(s^3)\right)\ve^-.\eeq

When the front crosses a rarefaction of the opposite family, its strength does
not change.

Finally,  when the  small wave impinges on a large shock at $P_1$ or at $P_2$,
we need to estimate  the relative size of the reflected wave front.

Calling $\ve^-$, $\ve^+$ the strengths of the front before and after interaction,
to leading order we have
\bel{reflec}\ve^+~=~\left(1- 2\tan\theta\right)\ve^-
\eeq
where $\theta$ is the angle between line segments $B_1 A_2$ and $A_1 A_2$
in Figure \ref{f:hyp59}, and $s$ is the strength of inner shocks $A_1C$ or $A_2C$.

\begin{figure}[htbp]
  \includegraphics[scale=0.5]{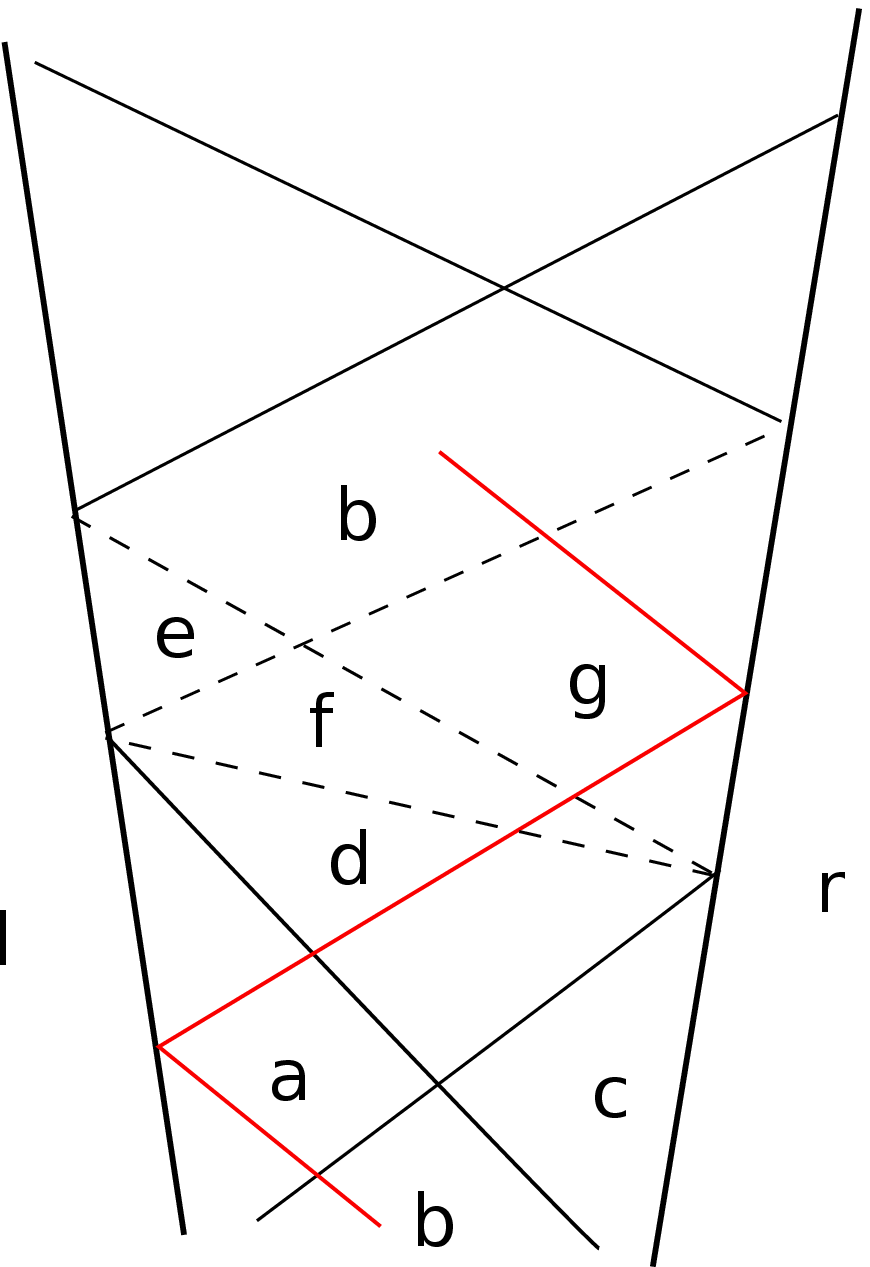}\qquad\quad 
  \includegraphics[scale=0.5]{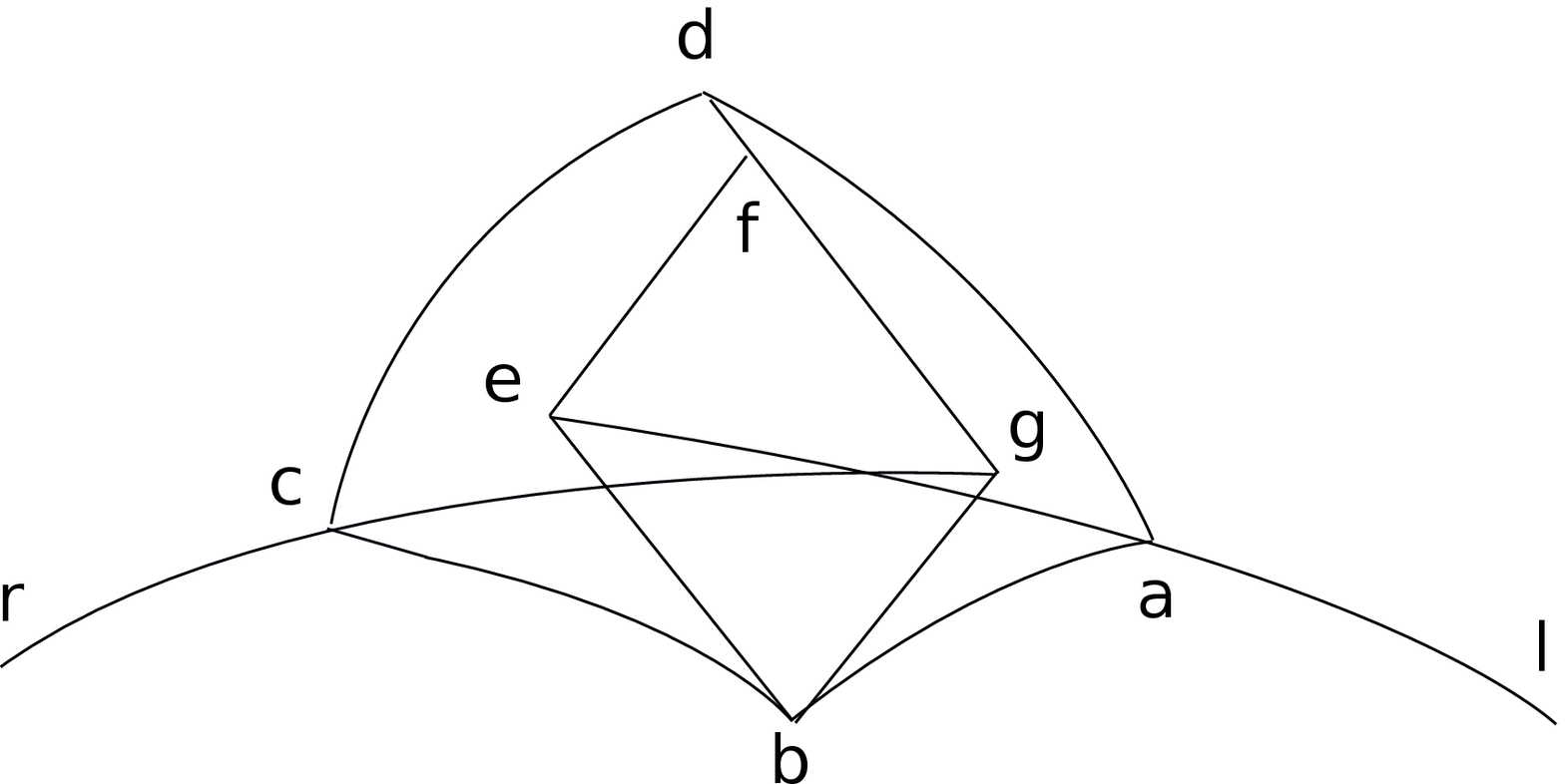}
    \caption{{\small Amplification of infinitesimal waves.
    }}
\label{f:hyp78}
\end{figure}

When the additional front reaches $B$, we want its size to be increased by a factor
$\kappa>1$. To achieve this goal, we need to perturb the
symmetric periodic pattern into an asymmetric periodic pattern
as shown in Figure \ref{f:hyp78}.

As in the figure, for simplicity we assume $h_c=1, u_c=0$, $u_d=r$, $u_b=s$.
Using the Rankine-Hugoniot  condition, we can calculate $h_b$ and $h_d$.

Indeed,
$u_b-u_c=s=\sqrt{-[p(v_b)-p(v_c)](v_b-v_c)}$.
\bel{s}
s=\sqrt{-[p'(v_b-v_c)+\frac12p''(v_b-v_c)^2+\frac16p^{(3)}(v_b-v_c)^3+\frac1{24}p^{(4)}(v_b-v_c)^4+\Ord ((v_b-v_c)^5)](v_b-v_c)}
\eeq
By expressing $v_b-v_c$ in powers of $s$, one obtains
$$v_b-v_c~=~1+\sqrt{-p'}s+\frac{p''}{4p'^2}s^2+\left(\frac5{32}\frac{p''^2}{(-p')^{7/2}}+\frac1{12}\frac{p^{(3)}}{(-p')^{5/2}}\right)s^3+\left(-\frac{p''^3}{8p'^5}+\frac{p''p^{(3)}}{8p'^4}-\frac{p^{(4)}}{48p'^3}\right)s^4+\Ord (s^5)
.$$

%
%
Considering \eqref{defh}, we have
\beq
h_b-1~=~\ds-s+\frac{1}{96}\frac{p''^2}{(-p')^3}s^3+
\left(\frac{p''^3}{64(-p')^{9/2}}+\frac{p''p^{(3)}}{96(-p')^{7/2}}\right)s^4+\Ord (s^5).
\eeq
In a similar way, we obtain
\beq
h_d-1~=~r-\frac{1}{96}\frac{p''^2}{(-p')^3}r^3+\left(\frac{p''^3}{64(-p')^{9/2}}+\frac{p''p^{(3)}}{96(-p')^{7/2}}\right)r^4+\Ord (r^5)
\eeq
Since $g$ is the intersection point of two rarefactions, we can calculate the coordinate of $g$ as
\beq
\bega{ll}
\ds u_g~=~(s+r)-\frac{1}{192}\frac{p''^2}{(-p')^3}(s^3+r^3)+\Ord (r^4,s^4),\cr\cr
\ds h_g~=~1+r-s+\frac{1}{192}\frac{p''^2}{(-p')^3}(s^3-r^3)+\frac12\left(\frac{p''^3}{64(-p')^{9/2}}+\frac{p''p^{(3)}}{96(-p')^{7/2}}\right)(r^4+s^4)+\Ord (r^5,s^5).
\enda
\eeq
Hence the slope of $cg$ is
\bel{theta0}
\ds\tan\theta~=~\frac{r-s+\frac{1}{192}\frac{p''^2}{(-p')^3}(s^3-r^3)+\frac12
\Big(\frac{p''^3}{64(-p')^{9/2}}+\frac{p''p^{(3)}}{96(-p')^{7/2}}\Big)(r^4+s^4)+\Ord (s^5,r^5)}{s+r-\frac{29}{192}\frac{p''^2}{(-p')^3}(s^3+r^3)}.
\eeq
For the left boundary shock, we can repeat above process. By assuming $u_a-u_b=\bar{s}, u_a-u_d=\bar{r}$, we can obtain the slope of $ae$, which is similar to \eqref{theta0},
\bel{theta}
\ds\tan\theta~=~\frac{\bar r-\bar s+\frac{1}{192}\frac{p''^2}{(-p')^3}(\bar s^3-\bar r^3)+\frac12(\frac{p''^3}{64(-p')^{9/2}}+\frac{p''p^{(3)}}{96(-p')^{7/2}})(\bar r^4+\bar s^4)+\Ord (\bar s^5,\bar r^5)}{\bar s+\bar r-\frac{29}{192}\frac{p''^2}{(-p')^3}(\bar s^3+\bar r^3)}.
\eeq
Here $r$ and $s$ are independent, so we can take different relations between $r(\bar r)$ and $s(\bar s)$ for the right and left shocks. From the figure, we expect $r\leq s$.

For the right boundary shock, we take
$r=s-(\frac{p''^3}{64(-p')^{9/2}}+\frac{p''p^{(3)}}{96(-p')^{7/2}})s^4$,
To leading order, the slope of the shock curve $cg$ is
$$
\tan\theta~=~o(s^3).
$$

We take $\bar r=\bar s-2J_2\bar s^4$, so the slope of $ae$ is
$$
\tan\theta~=~(-J_2+\frac12 (-p')^{-9/2}J_2)\bar s^3+o(\bar s^3).
$$

Since $r-s=\bar r-\bar s$, the relation between $s$ and $\bar s$ is
\beq
\frac1{6(-p')^{9/2}}J_2s^4=2J_2\bar s^4.
\eeq

After one complete set of interaction, the strength of the small wave located at $B$ is
$$
\ve_B~=~(1+J_2s^3+o(s^3))^2\left(1-2J_2\bar s^3+ (-p')^{-9/2}J_2\bar s^3+o(s^3)\right)(1+o(s^3))\ve_A~=~\left(1+ Xs^3+o(s^3)\right)\ve_A\,,
$$
where $$X=((-p')^{-9/2}-2)J_2\left(\frac{1}{12(-p')^{9/2}}\right)^{3/4}+2J_2>0.$$

The small wave has been amplified by a factor $1+ Xs^3+o(s^3)$.

By construction, after each period each pair of small compression-rarefaction wavefronts
is enlarged by a factor $\geq\lambda>1$.   When a pair grows to size $>\ve$, we can
perform a partial cancellation so that its size remains $\in [\ve/2,~\ve]$. After this manipulation, we can restrict the specific volume $v$ to be in the interval $(v_L,v_U)$. So the condition \eqref{blowup} always holds in the construction.

Since the total number of small wave-fronts is infinite,
after several periods a larger and larger number of  pairs (compression + rarefaction)
reaches size $>\ve/2$.   Hence, as $t\to\infty$,
 the total variation of this approximate solution grows without bounds.

%

\section*{Appendix}

Some detailed calculations about the slope of shock curves are given below.

As in the figure, for simplicity we assume $h_c=1, u_c=0$, $u_d=r$, $u_b=s$. By R-H condition,
$u_b-u_c=s=\sqrt{-[p(v_b)-p(v_c)](v_b-v_c)}$. By doing Taylor expansion, we have
\beq
s=\sqrt{-[p'(v_b-v_c)+\frac12p''(v_b-v_c)^2+\frac16p^{(3)}(v_b-v_c)^3+\frac1{24}p^{(4)}(v_b-v_c)^4+\Ord ((v_b-v_c)^5)](v_b-v_c)}
\eeq
We want to express $v_b-v_c$ in powers of $s$.
 Assume $v_b-v_c=1+As+Bs^2+Cs^3+Ds^4+\Ord (s^5)$, then compare the coefficients of $s^2, s^3, s^4, s^5, s^6$ to determine constants $A,B,C,D$.
\beq
\bega{ll}
s^2=-[p'(As+Bs^2+Cs^3+Ds^4+\Ord (s^5))+\frac12p''(As+Bs^2+Cs^3+Ds^4+\Ord (s^5))^2\cr\cr
+\frac16p^{(3)}(As+Bs^2+Cs^3+Ds^4+\Ord (s^5))^3+\frac1{24}p^{(4)}(As+Bs^2+Cs^3+Ds^4+\Ord (s^5))^4+\Ord (s^5)]\cr\cr
\cdot(As+Bs^2+Cs^3+Ds^4+\Ord (s^5))
\enda
\eeq

Coefficient for $s^2$: $$1=-A^2p',$$ $$A=\pm(-p')^{-1/2}.$$

Coefficient for $s^3$:$$0=-A(Bp'+\frac12A^2p'')-BAp',$$ $$B=-\frac{p''}{4p'}A^2=\frac{p''}{4p'^2}.$$

Coefficient for $s^4$:$$\bega{ll} 0&=-A(Cp'+p''AB+\frac16 p^{(3)}A^3)-B(Bp'+\frac12A^2p'')-CAp'\cr\cr&=-2ACp'-p''A^2B-\frac16p^{(3)}A^4-B^2p'-\frac12A^2Bp''.\enda$$
$$
2ACp'=\frac{p''^2}{4p'^3}-\frac{p^{(3)}}{6p'^2}-\frac{p''^2}{16p'^3}+\frac{p''^2}{8p'^3}=\frac{5p''^2}{16p'^3}-\frac{p^{(3)}}{6p'^2}
$$
$$C=\pm \frac5{32}\frac{p''^2}{(-p')^{7/2}}\pm\frac1{12}\frac{p^{(3)}}{(-p')^{5/2}}.$$

Coefficient for $s^5$:
$$
\bega{lll}
0&=&-A(Dp'+\frac12p''B^2+p''AC+\frac16p^{(3)}3A^2B+\frac1{24}p^{(4)}A^4)\cr\cr
&&-B(Cp'+p''AB+\frac16p^{(3)}A^3)\cr\cr
&&-C(Bp'+\frac12A^2p'')\cr\cr
&&-DAp'
\enda
$$
$$2ADp'=-2BCp'-\frac32p''AB^2-\frac32p''A^2C-\frac23p^{(3)}A^3B-\frac1{24}p^{(4)}A^5,
$$
$$
\bega{ll}
D&\ds=-\frac{BC}{A}-\frac{3}4\frac{p''}{p'}B^2-\frac34\frac{p''}{p'}AC-\frac13\frac{p^{(3)}}{p'}A^2B-\frac1{48}\frac{p^{(4)}}{p'}A^4\cr\cr
&\ds=-\frac{p''^3}{8p'^5}+\frac{p''p^{(3)}}{8p'^4}-\frac{p^{(4)}}{48p'^3}.
\enda
$$
\endproof

So by definition of $h_b$ and values of $A,B,C,D$, we have
\beq
\small{
\bega{lll}
h_b-1&=&\ds\int_{v_b}^{v_c}\sqrt{-p'}dV=-\int_{v_c}^{v_c+\Delta}\sqrt{-p'}dv\cr\cr
&=&-\ds\int_{v_c}^{v_c+\Delta}\sqrt{-p'}+\frac12(-p')^{-1/2}(-p'')(v-v_c)+\frac12[-\frac14(-p')^{-3/2}p''^2-\frac12(-p')^{-1/2}p^{(3)}](v-v_c)^2 \cr\cr
&&\ds+\frac16[-\frac38(-p')^{-5/2}p''^3-\frac12(-p')^{-3/2}p''p^{(3)}-\frac14(-p')^{-3/2}p''p^{(3)}-\frac12(-p')^{-1/2}p^{(4)}](v-v_c)^3
dv\cr\cr
&=&\ds-\sqrt{-p'}\Delta-\frac12(-p')^{-1/2}(-p'')\frac12\Delta^2-\frac16[-\frac14(-p')^{-3/2}p''^2-\frac12(-p')^{-1/2}p^{(3)}]\Delta^3\cr\cr
&&\ds -\frac1{24}[-\frac38(-p')^{-5/2}p''^3-\frac34(-p')^{-3/2}p''p^{(3)}-\frac12(-p')^{-1/2}p^{(4)}]\Delta^4\cr\cr
&=&\ds-(-p')^{1/2}As+\left(-(-p')^{1/2}B+\frac14(-p')^{-1/2}p''A^2\right)s^2\cr\cr
&&\ds+\left(-(-p')^{1/2}C+\frac12(-p')^{-1/2}p''AB-\frac16[-\frac14(-p')^{-3/2}p''^2-\frac12(-p')^{-1/2}p^{(3)}](-p')^{-3/2}\right)s^3\cr\cr
&&\ds+\left(-(-p')^{1/2}D-\frac12(-p')^{-1/2}(-p'')\frac12(B^2+2AC+\frac16[\frac14(-p')^{-3/2}p''^2
+\frac12(-p')^{-1/2}p^{(3)}]3A^2B\right.\cr\cr
&&\left.\ds-\frac1{24}[-\frac38(-p')^{-5/2}p''^3-\frac34(-p')^{-3/2}p''p^{(3)}-\frac12(-p')^{-1/2}p^{(4)}]A^4\right)s^4\cr\cr
&=&\ds-s+\frac{1}{96}\frac{p''^2}{(-p')^3}s^3+(\frac{p''^3}{64(-p')^{9/2}}+\frac{p''p^{(3)}}{96(-p')^{7/2}})s^4
\enda
}
\eeq

Since the rarefaction curves $dg$ and $bg$ are perpendicular in the Figure \ref{f:hyp78}, we can solve the following system for $(u_g,h_g)$.
\begin{equation}
\left\{\bega{l}
x-s=y-(1-s+\frac{1}{96}\frac{p''^2}{(-p')^3}s^3+(\frac{p''^3}{64(-p')^{9/2}}+\frac{p''p^{(3)}}{96(-p')^{7/2}})s^4+\Ord (s^5))\cr\cr
x-r=(1+r-\frac{1}{96}\frac{p''^2}{(-p')^3}r^3+(\frac{p''^3}{64(-p')^{9/2}}+\frac{p''p^{(3)}}{96(-p')^{7/2}})r^4+\Ord (r^5))-y
\enda\right.
\end{equation}
\bel{u_g}
\bega{ll}
\ds u_g=(s+r)-\frac{1}{192}\frac{p''^2}{(-p')^3}(s^3+r^3)+\Ord (r^4,s^4),\cr\cr
\ds h_g=1+r-s+\frac{1}{192}\frac{p''^2}{(-p')^3}(s^3-r^3)+\frac12\left(\frac{p''^3}{64(-p')^{9/2}}+\frac{p''p^{(3)}}{96(-p')^{7/2}}\right)(r^4+s^4)+\Ord (r^5,s^5).
\enda
\eeq
So the slope of the shock curve $rg$ at $g$ is
\bel{theta3}
\ds\tan\theta=\frac{r-s+\frac{1}{192}\frac{p''^2}{(-p')^3}(s^3-r^3)+\frac12(\frac{p''^3}{64(-p')^{9/2}}+\frac{p''p^{(3)}}{96(-p')^{7/2}})(r^4+s^4)+\Ord (s^5,r^5)}{s+r-\frac{29}{192}\frac{p''^2}{(-p')^3}(s^3+r^3)}
\eeq

\vskip 1em
{\bf Acknowledgments.} The research of the first author was partially supported
by NSF, with grant  DMS-1411786: ``Hyperbolic Conservation Laws and Applications".
The fourth author is supported in part
by National Natural Science Foundation of China under grant 11231006, Natural Science Foundation of Shanghai under grant 14ZR1423100 and  China Scholarship Council.
\vskip 1em

\end{document}